\documentclass[12pt,leqno]{amsart}
\usepackage{amsmath,amssymb,amsfonts}
\usepackage{eucal}

\newcommand \comment[1]{}			
\newcommand \dateadded[1]{\comment{[Date added: #1.]}}
\newcommand \mylabel[1]{\label{#1}\comment{{\rm \{#1\} }}}
\newcommand \myref[1]{\ref{#1}\comment{{\{#1\}}}}

\renewcommand{\phi}{\varphi}                 
\renewcommand{\epsilon}{\varepsilon}
\newcommand\eset{\varnothing}

\newtheorem{lem}{Lemma}[section]
\newtheorem{cor}[lem]{Corollary}

\newtheorem{thm}[lem]{Theorem}

\numberwithin{equation}{section}

\newcommand\inv{^{-1}}
\newcommand\textb{\text{\rm b}}
\newcommand\chib{\chi^{\textb}}
\newcommand\Lat{\operatorname{Lat}}
\newcommand\Latb{\operatorname{\Lat^{\textb}}}

\newcommand\cA{\mathcal{A}}
\newcommand\cB{\mathcal{B}}
\newcommand\cL{\mathcal{L}}
\newcommand\cS{\mathcal{S}}

\newcommand\bbR{\mathbb{R}}
\newcommand\bbZ{\mathbb{Z}}

\newcommand \fA{\mathfrak A}

\newcommand \bL{\mathbf L}

\newcommand\M{m}
\newcommand\chiZ{\chi^\bbZ}
\newcommand\pZ{p^\bbZ}
\renewcommand\mod{\operatorname{mod}}
\newcommand\chimod{\chi^{\mod}}
\newcommand\G{\Gamma}
\renewcommand\l{\lambda}
\renewcommand\L{\Lambda}
\newcommand\LOOP{\operatorname{lp}}

\begin{document}

\begin{center}

\large
\textsc{Lattice Point Counts for the Shi Arrangement\\ and Other Affinographic Hyperplane Arrangements}
\normalsize
\vskip20pt

{David Forge\footnote{The research of the first author was performed while visiting the State University of New York at Binghamton.}\\
Laboratoire de recherche en informatique UMR 8623\\
B\^at.\ 490, Universit\'e Paris-Sud\\
91405 Orsay Cedex, France\\
E-mail: {\tt forge@lri.fr}}\\[10pt]

and\\[10pt] 

{Thomas Zaslavsky\\
Department of Mathematical Sciences \\
State University of New York at Binghamton\\
Binghamton, NY 13902-6000, U.S.A.\\
E-mail: {\tt zaslav@math.binghamton.edu}}\\[20pt]

{Version of \today.}\\[20pt]

\end{center}

\small{\sc Abstract.}
 Hyperplanes of the form $x_j = x_i + c$ are called \emph{affinographic}. 
 For an affinographic hyperplane arrangement in $\bbR^n$, such as the Shi arrangement, we study the function $f(\M)$ that counts integral points in $[1,\M]^n$ that do not lie in any hyperplane of the arrangement.  We show that $f(\M)$ is a piecewise polynomial function of positive integers $\M$, composed of terms that appear gradually as $\M$ increases.  Our approach is to convert the problem to one of counting integral proper colorations of a rooted integral gain graph.
 
 An application is to interval coloring in which the interval of available colors for vertex $v_i$ has the form $[h_i+1,\M]$.

 A related problem takes colors modulo $\M$; the number of proper modular colorations is a different piecewise polynomial that for large $\M$ becomes the characteristic polynomial of the arrangement (by which means Athanasiadis previously obtained that polynomial).  We also study this function for all positive moduli.
\dateadded{15 July 2005}
\normalsize\\

\emph{Mathematics Subject Classifications (2000)}:
{\emph{Primary} 05C22, 52C35; \emph{Secondary} 05C15.}\\

\emph{Key words and phrases}:
{Integral gain graph, modular gain graph, proper coloring, interval graph coloring, chromatic function, affinographic hyperplane arrangement, deformation of Coxeter arrangement, Shi arrangement, Linial arrangement}\\


\vskip 20pt

 \section{Integral affinographic hyperplane arrangements}\mylabel{intro}

The Shi arrangement of hyperplanes is the set $\cS_n$ of affine hyperplanes in $\bbR^n$ that consists of the hyperplanes with equations $x_j=x_i$ and $x_j=x_i+1$ for all pairs $i<j$ in $\{1,2,\dots,n\}$.  
A \emph{region} of $\cS_n$ is a component of the complement $\bbR^n \setminus \big(\bigcup\cS_n\big)$.  
The number of regions formed by the Shi arrangement has some interest that does not directly concern us (see \cite{ShiMono}), but which led to a series of investigations, first to find that number \cite{Shi}, then to find the more refined invariant called the characteristic polynomial, $p(\l)$, from which the number of regions is a simple deduction, and then to characteristic polynomials of analogous arrangements.
For instance, for the Shi arrangement the characteristic polynomial is $\l(\l-n)^{n-1}$ \cite{Headley}; and one of the important analogs is the Linial arrangement, $\{x_j=x_i+1 : i<j\}$.  

 There is an interpretation of $p(\l)$ in terms of graph theory.  An \emph{integral gain graph} is a graph with edges labelled invertibly by integers.  (Section \myref{gg} has the precise definition.)  For each arrangement of integral affinographic hyperplanes there is an associated integral gain graph.  For any sufficiently large positive integer $\M$, $p(\M)$ is the number of ways to properly color this graph by colors in the cyclic group $\bbZ_\M$.  (In saying this we follow both the implicit content of Athanasiadis \cite{Athan}, as reformulated in Section \myref{modular} below, and also \cite{PDS}.)  
That is analogous to ordinary graph coloring, where the chromatic polynomial $\chi_\G(\l)$ of a graph $\G$, evaluated at any positive integer $\M$, equals the number of proper colorations with values in $\bbZ_\M$.  

Viewing ordinary graph coloring geometrically leads to two other interpretations of the chromatic polynomial.  First, to $\G$ there corresponds a hyperplane arrangement that consists of all the hyperplanes $x_i=x_j$ for which there is an edge $e_{ij}$ in $\G$; then $\chi_\G(\l)$ equals the characteristic polynomial of the arrangement, $p(\l)$.  Second, for every $\M > 0$, $p(\M)=\chi_\G(\M)$ is also the number of integer points that are in the hypercube $[1,\M]^n \subset \bbR^n$ but not in any hyperplane.  

For the Shi and Linial arrangements this last interpretation is invalid.  Instead, the number of integer points is a new function, $\pZ(\M)$, that does not agree with $p(\M)$. 
 This new function is our main topic.  We investigate it for a wide class of hyperplane arrangements, known variously as integral affinographic arrangements or as integral deformations of the complete-graph arrangement $A_{n-1}^*$.  
We treat it through the interpretation of an integral affinographic arrangement as an integral gain graph.  (Thus our method resembles that of Athanasiadis.) 
We compute the \emph{integral chromatic function} $\chiZ(\M)$, defined as the number of proper colorations using the color set $\{1,2,\dots,\M\}$ where $\M\geq0$, in two theoretical ways: first, by M\"obius inversion over the semilattice of balanced flats of a matroid on the gain graph, and second, by a deletion-contraction formula.  We also develop examples, including the extended Shi and Linial arrangements, where we count colorations combinatorially.  

Our central idea is to transfer the problem to \emph{rooted} gain graphs, which makes it easy to keep track of the various graph transformations necessary to compute the integral chromatic function.  The main theorem is that this function is a sum of terms, added in gradually as $\M$ increases, each term being a polynomial with integral roots that is nonnegative for those values of $\M$ for which it is included in the evaluation of $\chiZ(\M)$.

Our method is adaptable to modular coloring of affinographic arrangements, thus giving a formula for the number of $\bbZ_\M$-colorations even for small moduli where Athanasiadis' general polynomial formula does not apply.  We do this in Section \myref{modular}.

It is also applicable to a special kind of list coloring of ordinary graphs.  Suppose to each vertex $v_i$ we assign an integer $h_i$ and require that $v_i$ be colored with a color in the interval $(h_i, \M] \subseteq \bbZ$.  The number of proper colorations is a function of $\M$ that is zero for large negative $\M$ and equals the chromatic polynomial for large positive $\M$.  In between, its behavior is like that of the integral chromatic function.  From the geometric viewpoint, we are counting the integral points $x$ in the parallelepiped $(h_1,\M] \times\cdots\times (h_n,\M]$ that avoid all the hyperplanes corresponding to the edges of the graph.  This situation has curious similarities to a theory of Noble and Welsh \cite{NW}.  See Section \myref{graph}.
\dateadded{2005 Aug 2}

\section{Background}\mylabel{background}

\subsection{Basic definitions}\mylabel{defs}

We adopt the notation $[n]:=\{1,2,\dots,n\}$ for a positive integer $n$; and $[0] := \eset$.

A graph $\G=(V,E)$ may have loops and multiple edges (but not the half and 
loose edges of \cite{BG}).  A \emph{link} is an edge that has two 
distinct endpoints; thus, our edges are links and loops.  A \emph{circle} 
is a connected subgraph of degree 2, or its edge set.\footnote{We eschew the common term ``cycle'' due to its several inconsistent meanings in graph theory.\dateadded{15 Feb 2006}}  
We may write a circle $C$ as a word $e_1e_2\cdots e_l$; this means the edges are numbered 
consecutively around $C$ and oriented in a consistent direction.

\subsection{Gain graphs}\mylabel{gg}

An \emph{abelian gain graph}, $\Phi=(\G,\phi,\fA)$, consists of a graph 
$\G$, an abelian group $\fA$ called the \emph{gain group}, and a 
\emph{gain function} $\phi : E \to \fA$ that is orientable, i.e., if $e$ 
denotes an edge oriented in one direction and $e\inv$ the same edge with 
the opposite orientation, then $\phi(e\inv)=-\phi(e)$.  The gain of a 
circle $C=e_1e_2\cdots e_l$ is $\phi(C) = 
\phi(e_1)+\phi(e_2)+\cdots+\phi(e_l)$.  This is not entirely well defined 
since it depends on the direction but it is well defined whether $C$ has 
zero or nonzero gain.  We define $\cB(\Phi) := \{ C : \phi(C)=0 \}$.  A circle in 
$\cB$, and more broadly any subgraph or edge set all of whose circles are 
in $\cB$, is called \emph{balanced}.

Our concern will be principally with \emph{integral gain graphs}, whose   
gain group is the additive group of integers, and to a lesser extent with 
\emph{modular gain graphs}, whose gain group is $\bbZ_\M$.

For $S\subseteq E$ we write $\Phi|S$ for the subgraph $(V,S)$ with gains as in $\Phi$, $\pi(S) := \pi(\Phi|S) :=$ the partition of $V$ whose blocks are the vertex sets of the components of $\Phi|S$, and $b(S) := b(\Phi|S) :=$ the number of balanced components of $\Phi|S$.  When $S$ is balanced, $b(S) = |\pi(S)|$.
 The notation $\LOOP(\Phi)$ means $\Phi$ with all links deleted.
 A convenient notation for an edge is $ge_{ij}$; this means the endpoints are $v_i$ and $v_j$ and the gain in the direction from $v_i$ to $v_j$ is $g$.  Generally, $e_{ij}$ denotes an edge whose endpoints are $v_i$ and $v_j$, oriented from $v_i$ to $v_j$.
 Since we have many parallel edges, sometimes it is convenient to indicate them all by the notation $\fA_{ij}e_{ij}$, that is, a single edge $e_{ij}$ labelled with the entire set $\fA_{ij}$ of gains of edges $ge_{ij}$.  However, one should keep in mind that this is merely notational shorthand and each separate edge $ge_{ij}$ continues to be a distinct object.

A \emph{group coloration} of $\Phi$ is a mapping $x: V \to \fA : v_i \mapsto x_i$.  It is \emph{proper} if, whenever there is an edge $ge_{ij}$, then $x_j \neq x_i+g$; in general the set of \emph{improper} edges of $x$ is 
$$I(x) := \{ ge_{ij} : x_j=x_i+g \}.$$
 If $\fA$ has finite order $\M$ we can consider the number of proper group colorations of $\Phi$.  By taking new gain groups that are supergroups of $\fA$ we obtain a function of the group; this function is a polynomial in the order of the group, called the \emph{balanced chromatic polynomial} of $\Phi$ and written $\chib_\Phi(\l)$; it has the formula
 \begin{equation}\mylabel{E:balchromatic}
 \chib_\Phi(\l) = \sum_{S\subseteq E : \text{ balanced}} 
(-1)^{|S|} \l^{b(S)} .
 \end{equation}
(Ordinary graphs can be regarded as the special case in which all gains are zero.  Then $\chib_\Phi$ equals the usual chromatic polynomial of $\G$.)
 Even if $\fA$ is infinite, as long as $\Phi$ is finite one can define $\chib_\Phi(\l)$ by Equation \eqref{E:balchromatic} or by a second algebraic expression for the balanced chromatic polynomial in terms of $\Latb \Phi$, the class of all closed, balanced edge sets, ordered by inclusion.  A balanced edge set $B$ is called \emph{closed} if any edge $ge_{ij}$ whose endpoints are joined by a path $P$ in $B$ with the same gain (that is, $P$ is open, possibly of length $0$, and has gain $g$ in the orientation from $v_i$ to $v_j$) is itself in $B$.  Then 
 \begin{equation}\mylabel{E:balchromaticmu}
 \chib_\Phi(\l) = \sum_{B\in\Latb\Phi} \mu(\eset,B) \l^{b(B)} ,
 \end{equation}
 where $\mu$ is the M\"obius function of $\Latb\Phi$.  We mention that $I(x)$ is closed and balanced.  
 (These facts are adaptations of results in \cite[Sections III.4 and III.5]{BG}.  Colorations in a supergroup are enumeratively equivalent to the zero-free $k$-colorations of \cite[Section III.4]{BG}, $k$ being the index of $\fA$ in the supergroup.  For the M\"obius function of a poset see \cite{FCT,EC1}.)

If $\Psi$ is a different gain graph with the same underlying graph $\G$ and the same list of balanced circles (that is, $\cB(\Psi) = \cB(\Phi)$), then \eqref{E:balchromatic} shows that $\chib_\Psi = \chib_\Phi$.  In particular, suppose $\fA=\bbZ$.  Then we can replace $\bbZ$ by a sufficiently large finite cyclic group $\bbZ_\M$.  See \cite[Section 11.4]{PDS} or Section \myref{modular} for applications of this idea.

\subsection{Switching and potentials}\mylabel{switching}

There is a transformation of gain graphs called \emph{switching} that does not change the balanced circles.  To switch $\Phi$ we take a \emph{switching function} $\eta: V \to \fA : v_i \mapsto \eta_i$ and replace $\phi$ by $\phi^\eta$, whose definition is
 $$
 \phi^\eta(e_{ij}) := \phi(e_{ij}) + \eta_j - \eta_i .
 $$
 We write $\Phi^\eta := (\G,\phi^\eta,\fA)$.  It is clear that switching does not change balance or the balanced chromatic polynomial.  

Suppose we have an edge set $S$ such that some switching of $\Phi$ converts all the gains on $S$ to zero; obviously then $S$ is balanced.  Conversely, if $S$ is balanced in $\Phi$ there is a switching function $\eta$ such that the gains on $S$ are all zero in $\Phi^\eta$ \cite[Section I.5]{BG}.  Indeed, we can specify $\eta$: if there is a path $P$ in $S$ from $v_i$ to $v_j$, then $\eta_j$ must equal $\eta_i+\phi(P)$.  This rule leaves one value of $\eta$ to be chosen arbitrarily in each component of $\Phi|S$, after which $\eta$ is fully determined.

The negative of a switching function for $S$ is called a \emph{potential} for $S$; that is, a potential for $S$ is a mapping $\theta : V \to \fA$ such that $\phi(e_{ij}) = \theta_j-\theta_i$ for any edge in $S$.  A group coloration $x$ is a potential for $I(x)$.

Switching acts on group colorations in such a way as to maintain propriety.  A coloration $x$ of $\Phi$ switches to $x^\eta$ defined by $x^\eta_i = x_i+\eta_i$.  Then $I(x^\eta)=I(x)$, where $I(x^\eta)$ is calculated in $\Phi^\eta$ and $I(x)$ is calculated in $\Phi$.

\subsection{Contraction}\mylabel{contraction}

\emph{Contraction} of $\Phi$ by a balanced edge set $S$ consists of two steps.  First, $\Phi$ is switched so that $S$ has all zero gains.  Then, each block $W\in\pi(S)$ is identified to a single vertex and $S$ is deleted.  
 The notation for $\Phi$ with $S$ contracted is $\Phi/S$.

\subsection{Affinographic arrangements and their gain graphs}\mylabel{affinogg}

An \emph{integral affinographic hyperplane arrangement} in $\bbR^n$ (we shall omit the words ``integral'' and ``hyperplane'') is a finite set $\cA$ of affine hyperplanes of the form $h_{ij}(g) : x_j-x_i=g$ for $i,j\in[n]$ and $g\in\bbZ$.  The \emph{intersection semilattice} of $\cA$ is the set $\cL(\cA)$ of all nonempty intersections of subsets of $\cA$ (including the intersection of no hyperplanes, which is $\bbR^n$); it is partially ordered by reverse inclusion.  It is a meet semilattice with $\hat0=\bbR^n$ (in fact, a geometric semilattice \cite{McNaff,WW}) and it has a $\hat1$ if and only if there is a point common to every hyperplane.  The \emph{characteristic polynomial} of $\cA$ is
 $$
 p_\cA(\l) := \sum_{W\in\cL(\cA)} \mu(\hat0,W) \l^{n-\dim W} ,
 $$
 where $\mu$ denotes the M\"obius function of $\cL(\cA)$. 
 One of the interesting properties of $p_\cA$ is that $(-1)^n p_\cA(-1)$ is the number of regions formed by $\cA$.

To $\cA$ there corresponds an integral gain graph $\Phi(\cA)$.  The vertex set is $V = \{v_1,v_2,\ldots,v_n\}$.  The gain group is $\bbZ$.  Corresponding to $h_{ij}(g) \in \cA$ there is an edge $ge_{ij}$.  Thus, the gain function is given by $\phi(ge_{ij})=g$.  By the definition of gains, the gain in the opposite direction is $-g$, which corresponds to the fact that $h_{ij}(g)$ and $h_{ji}(-g)$ are the same hyperplane.  Thus there are one vertex for each coordinate and one edge for each hyperplane.

The significance of this correspondence is, first, that the characteristic polynomial $p_\cA$ equals the balanced chromatic polynomial of $\Phi$ (by \cite[Theorem III.5.2 and Corollary IV.4.5]{BG}); thus we can compute the former by coloring $\Phi$.  But for us the more important fact is that the lattice points we wish to count are exactly the same as the proper integral colorations of $\Phi(\cA)$.

\section{Integral coloring}\mylabel{coloring}

An \emph{(integral) $\M$-coloration} of an integral gain graph $\Phi$ is a function 
$x : V \to [\M] \subseteq \bbZ$; that is, it is a group coloration with a restricted range.

The gain group $\bbZ$ being an ordered group, it is possible to single out a canonical switching function $\eta$ for each balanced edge set $B$ by choosing $\eta$ so that it switches $B$ to have all zero gains and in each block $W \in \pi(B)$ its minimum value is $0$.  We call this $\eta$ the \emph{top-vertex switching function} for $B$.  
The maximum of $\eta$ in $W$, written $h_B(W)$, is called the \emph{height} of $W$ (or of the corresponding component of $B$).  A simpler description of $h_B(W)$ is that it equals the maximum gain of a path in that component of $B$ of which $W$ is the vertex set.  
 A vertex in $W$ with $\eta(v)=0$ is called a \emph{top vertex} (because it has maximum potential).  From now on, in contracting an integral gain graph by a balanced edge set we switch by the top-vertex switching function.  We call this \emph{top-vertex switching}.
 In top-vertex switching and contraction on $\Phi$ one may think of the contracted graph $\Phi/B$ as having vertex set obtained by retaining one top vertex $v_i$ from each $W\in\pi(B)$, the other vertices of $W$ being contracted into $v_i$.  
 The effect of top-vertex switching on a coloration is to change $x_j$ to $x_j+\eta_j$.

In the statement of the theorem we need the positive part of a real number $r$, which is $r^+ := \max(r,0)$.

\begin{thm}\mylabel{T:integralexpansion}
For $\M\in\bbZ$, 
 \begin{equation}\mylabel{E:integralexpansion}
 \chiZ_\Phi(\M) = \sum_{B\in\Latb\Phi} \mu(\eset,B) \prod_{W\in\pi(B)} 
\big[\M-h_B(W)\big]^+ ,
 \end{equation}
where $\mu$ is the M\"obius function of $\Latb\Phi$.
 \end{thm}

The proof will come shortly, but first we explore the meaning of the expression for $\chiZ_\Phi$.  Each $h_B(W)$ is a nonnegative integer and, if $W$ is a singleton $\{w\}$, then $h_B(W)=0$.  We thereby obtain a general description of the form of $\chiZ_\Phi$.
 
\begin{cor}\mylabel{C:integralexpansion}
 Let $\Phi$ be a nonempty integral gain graph of order $n$.  There exist a positive integer $k$, positive integers $d_1 = n > d_2, \ldots, d_k$, polynomials $p_j(\l) = (\l-r_{j1}) (\l-r_{j2}) \cdots (\l-r_{jd_j})$ of degree $d_j$ with integral roots $r_{j1} \geq r_{j2} \geq \cdots \geq r_{jd_j} \geq 0$ and with $p_1(\l) = \l^n$, and positive integers $\mu_1=1,\mu_2,\ldots,\mu_k$, such that 
 $$
 \chiZ_\Phi(\M) = \M^n + \sum_{2 \leq j \leq k:\, r_{j1} \leq \M} (-1)^{n-d_j}  \mu_j p_j(\M)
  \quad\text{ for } m\geq0 .
 $$
 \end{cor}
 
\begin{proof}
 Each $p_j(\l)$ corresponds to a balanced flat $B$; $\mu_j$ is the unsigned M\"obius function $|\mu(\eset,B)|$.  The degree $d_j = |\pi(B)|$, and the integers $r_{ji}$ are the numbers $h_B(W)$ for $W \in \pi(B)$.  We choose $p_1(\l)$ to correspond to $B = \eset$, so $p_1(\l) = \l^n$.  The rest is clear from Theorem \myref{T:integralexpansion}, the facts that every interval in $\Latb\Phi$ is a geometric lattice and that $B$ has rank $n - |\pi(B)|$, and the consequence by Rota's theorem that $(-1)^{n-d_j} \mu(\eset,B) > 0$.
 \end{proof}
 
Thus, the domain $[0,\infty)$ breaks up into intervals, on each of which $\chiZ_\Phi$ is a different monic polynomial of degree $n$.  Progressing up the domain, on each higher interval more terms $(-1)^{n-d_j} \mu_j p_j(\M)$ are added to $\chiZ_\Phi$.  Once term $j$ is added, which happens at $\M = r_{j1} = \max r_{ji}$, it remains as $\M$ continues to increase.  Each $r_{ji}$, for fixed $j$, is the largest gain of a path in a different component of the flat $B$ corresponding to term $j$; thus, in particular, the last term appears just when $\M$ reaches the maximum gain of a path in $\Phi$.  Each $\mu_j p_j$ is a polynomial whose values are positive for all $\M > r_{j1}$.  (The foregoing description translates in \cite{QRS} into an interesting aspect of the solution of a recreational chess problem.)
 \dateadded{Corollary etc.: 2005 Aug 2--11}

To prove the theorem we need our key innovation.  A \emph{rooted integral gain graph} is an integral gain graph with a distinguished \emph{root} vertex $v_0$, whose incident edges are called \emph{root edges}.  
The root vertex and root edges are subject to three rules.  
First, for each nonroot vertex $v_i$, the gains of root edges $e_{0i}$ form a lower interval of integers, $(-\infty,h_i] \subset \bbZ$ for some integer $h_i$.  Second, a root edge cannot be deleted or contracted.  Finally, in integral coloring the root vertex is always colored $0$.  The effect is that the color of a nonroot vertex $v_i$ must be greater than $h_i$.

Given an unrooted integral gain graph $\Phi$, its \emph{rooting} $\Phi_0$ is $\Phi$ with a new root $v_0$ adjoined, together with edges joining $v_0$ to all other vertices $v_i$, carrying all possible nonpositive gains in the direction $v_0v_i$.  In set shorthand, all the edges from the root to $v_i$ in the rooting of $\Phi$ are indicated by the notation $(-\infty,0]e_{0i}$.  
Rooting imposes a lower bound of $1$ on the color of every vertex in $\Phi$.

Not every rooted integral gain graph is a rooting of a gain graph, but, with the appropriate definitions, the essence of the theorem is valid for any rooted integral gain graph.  
A \emph{rooted (integral) $\M$-coloration} of the rooted gain graph $\Psi$ is any mapping $x: V(\Psi) \to (-\infty,\M]$ that obeys the coloring rule $x_0=0$ and for which every root edge is proper.  
The definitions of propriety of an edge and of a coloration are the same for rooted as for unrooted integral gain graphs, so while every rooted coloration of $\Psi$ is proper on root edges, it is the \emph{proper} rooted colorations that are proper on all edges.
(In particular, an $\M$-coloration of an unrooted integral gain graph $\Phi$ is the same, with the addition of the obligatory color $0$ at the root, as a rooted $\M$-coloration of its rooting, and the proper colorations of $\Phi$ and its rooting are also the same; thus, $\chiZ_{\Phi_0}(\M) = \chiZ_\Phi(\M)$.)  
The \emph{improper edge set} $I(x)$ of a rooted $\M$-coloration $x$ of $\Psi$ is the same as that of $x$ restricted to nonroot vertices.  This set is closed and balanced in $\Psi\setminus v_0$ as well as in $\Psi$.  
For a balanced set $B$ in $\Psi\setminus v_0$, we define $\pi_0(B)$ to be the corresponding partition of $V(\Psi)$ and $\pi(B)$ to be that of $V(\Psi)\setminus v_0$.  
A top vertex of $B$ is (as before) a vertex $v_i$ such that no path in $B$ that ends at $v_i$ has negative gain.  The top-vertex switching function $\eta$ is given by $\eta_j = \phi(P)$, if $P$ is a path in $B$ from $v_j$ to a top vertex; hence
$$
\eta_j = \max\{\phi(P) : P \text{ is a path in $B$ starting at } v_j\},
$$
and switching $\Psi$ to $\Psi^\eta$ yields $h_j^\eta = h_j + \eta_j$.  From these values we obtain a general definition of $h_B$ by setting
$$
h_B(W) := \max_{v_j\in Q} (h_j+\eta_j)
$$
for each set $W \in \pi(B)$.  
Theorem \myref{T:integralexpansion} generalizes to

\begin{thm}\mylabel{L:rootedintegralexpansion}
For any rooted integral gain graph $\Psi$ and integer $\M$, 
 \begin{equation}
 \chiZ_\Psi(\M) = \sum_{B\in L} \mu_L(\eset,B) \prod_{W\in\pi(B)}
\big[\M-h_B(W)\big]^+ ,
 \end{equation}
where $L:=\Latb(\Psi\setminus v_0)$.
\end{thm}

\begin{proof}
 Consider all rooted $\M$-colorations of $\Psi$, proper or not.  It follows 
by M\"obius inversion, as in \cite[p.\ 362]{FCT} or \cite[Theorem 2.4]{SGC}, that
 $$
 \chiZ_\Psi(\M) = \sum_{B\in L} \mu_L(\eset,B) f(B) ,
 $$
 where $f(B)$ is the number of rooted colorations $x$ such that $I(x)=B$. 

 We now show by a bijection that $f(B)$ is the number of proper $\M$-colorations of $\Psi/B$.  Let $\eta$ be the top-vertex switching function for $B$; in particular, $\eta_0=0$.  Switching an $\M$-coloration $x$ of $\Psi$ gives an $\M$-coloration of $\Psi^\eta$ (and conversely), because the color $x_j$ of a vertex with $\eta_j\neq0$ is changed to the color $x_i = x_j+\eta_j$ of a top vertex of the set $W\in\pi_0(I(x))$ that contains $v_j$, and $\M \geq x_i \geq x_j$.  
Also, the set of improper edges remains the same.  Therefore, when we contract $\Psi$ by $I(x)$ we get an $\M$-coloration with no improper edges.

Conversely, for any $B\in L$, let $\eta$ be the top-vertex switching function.  A proper rooted $\M$-coloration $y$ of $\Psi/B$ pulls back to an $\M$-coloration of $\Psi^\eta$ by $x_i = y_W$ where $v_i\in W \in \pi_0(B)$.  Then switching back to $\Psi$ we have an $\M$-coloration $x^{-\eta}$ of $\Psi$ whose improper edge set is $B$.

Since it is clear that these correspondences are inverse to each other, the bijection is proved.
 \end{proof}

\begin{proof}[Proof of Theorem \myref{T:integralexpansion}.]  
We apply the preceding theorem to the rooting of $\Phi$.
 \end{proof}

Let us write the formal integral combination
 \begin{equation}\mylabel{E:coloringsum}
 \bL_0(\Psi) := \sum_{B \in L} \mu_L(\eset,B) \LOOP(\Psi/B) .
 \end{equation}
 Then we can interpret the last theorem as a rule for evaluating $\bL_0(\Psi)$, and the theorem as an evaluation of $\bL_0(\Phi_0)$.  To state the rule, let us define the product of rooted graphs as their amalgamation at the root vertex.  Then the evaluation rule is that a rooted integral gain graph with a single nonroot vertex $v$ and root-edge gain set $(-\infty,h]$ evaluates to $0$ if $v$ supports a loop with gain $0$ and to $(\M-h)^+$ otherwise.  This way of thinking suggests a theory of Tutte invariants, which we begin to develop in \cite{PTRI}.

There is another way to express the number of proper $\M$-colorations: by deletion and contraction; but it applies only in the general setting of rooted gain graphs.

\begin{thm}\mylabel{L:rooteddandc}
Let $\Psi$ be a rooted integral gain graph and $e$ a nonroot link.  Then
 $$
 \chiZ_\Psi(\M) = \chiZ_{\Psi\setminus e}(\M) - \chiZ_{\Psi/e}(\M) .
 $$
\end{thm}

\begin{proof}
 The standard method works: we consider those proper $\M$-colorations of 
$\Psi\setminus e$ for which $e$, when restored to $\Psi$, is a proper edge, 
and those for which it is improper.  The former are proper 
$\M$-colorations of $\Psi$ and the latter correspond to proper 
$\M$-colorations of $\Psi/e$ because of the way contraction affects 
colorations, as discussed in the proof of Theorem 
\myref{L:rootedintegralexpansion}.
 \end{proof}

There is no corresponding result for $\M$-colorations of an unrooted integral gain graph.  In general, 
 $
 \chiZ_\Phi(\M) \ne \chiZ_{\Phi\setminus e}(\M) - \chiZ_{\Phi/e}(\M) ,
 $
 because the lower bound of $1$ on the color of a vertex in $\Phi$ and $\Phi\setminus e$ changes at the contracted vertex of $\Phi/e$ in the course of top-vertex switching.

The two kinds of formula we have given are related through broken balanced circles.  Given a linear ordering of the edges, a \emph{broken balanced circle} is a balanced circle without its last edge.  
(This is a special kind of broken circuit; we are adapting the theory of no-broken-circuit sets to geometric semilattices, specialized to the case of graphic lift matroids \cite[Section II.3]{BG}.)  
When using deletion and contraction to compute the number of proper $\M$-colorations, one linearly orders the edge set and, in sequence from first to last (except for edges that have become loops), contracts and deletes edges in every possible way.  If in this process a balanced circle is contracted to a loop, the resulting graph will have no proper colorations and will contribute $0$ to the total number.  The only way to avoid this and get a positive contribution is for the set $F$ of contracted edges to contain no broken balanced circle.  
$F$ is therefore a forest that contains no broken balanced circles.  One may conclude that the coefficient $\mu(\eset,B)$ is, up to sign, the number of forests $F \subseteq B$, with partition $\pi(F)=\pi(B)$, which contain no broken balanced circle.  
(We omit details, which are as in the standard broken-circuit theory.)

\section{The geometry of integral coloring}\mylabel{geom}

By thinking of an integral coloration $x$ as a point of the integral lattice $\bbZ^n$ in $\bbR^n$ we obtain the main theorem.

\begin{thm}\mylabel{T:latticepoints}
Let $\cA$ be an integral affinographic hyperplane arrangement in $\bbR^n$ and let $\M \in \bbZ$ be nonnegative.  The number of integer points in $[\M]^n$ that are contained in none of the hyperplanes of $\cA$ equals $\chiZ_{\Phi(\cA)}(\M)$ in Theorem \myref{T:integralexpansion}.
\end{thm}

With geometry we can do more: we can separately interpret each term of Equation \eqref{E:integralexpansion}.  
In $\bL_0(\Phi)$ each term is an integral gain graph $\L=\LOOP(\Phi_0/B)$ with a positive or negative  weight $\mu(\eset,B)$.  Taking the viewpoint that the vertices of $\L$ are the top vertices $v_i$ of the components of $\Phi|B$, $\L$ has links $e_{0i}$ with gain sets $(-\infty,h_i]$; thus the color of $v_i$ is 
restricted by the bound $x_i>h_i\geq0$.  
Furthermore, if any other vertex $v_j$ was contracted into $v_i$, it was contracted along a path $e_{jj_1}e_{j_1j_2}\cdots e_{j_ki}$ with total gain $g_{ji}$, say, that corresponds to the equation $x_i=x_j+g_{ji}$.  
This leads us to define for each $\L$ the relatively open cone
 \begin{align*}
 C(\L) := \{ x\in\bbR^n :\ & x_i > h_i \text{ for } v_i\in V(\L) \text{ and } \\
    & x_j=x_i-g_{ji} \text{ for all other vertices } v_j \}.
 \end{align*}
 We assign to $C(\L)$ the weight of $\L$, that is, $\mu(\eset,B)$.

 \begin{thm}\mylabel{T:conecolors}
 The proper integral colorations of $\Phi$ are the integral points $x$ in the positive orthant $\bbR_{>0}^n$ whose total weight, summed over all cones $C(\L)$ that contain $x$, is nonzero; and each of these points has total weight equal to $1$.
 \end{thm}

Counting only points $x \in [\M]^n$, we recover Theorem \myref{T:latticepoints}.

 \begin{proof}
 Consider a positive integral point $x\in\bbR^n$.  As a coloration of $\Phi$ it has an improper edge set $I(x)$.  
 The cones $C(\LOOP(\Phi_0/B))$ to whose affine span $x$ belongs are the ones for which $B \leq I(x)$ in $\Latb\Phi$.  Then $x$, being a potential for $I(x)$, is also a potential for $B$.  If $x$ is in the cone of every $B \leq I(x)$, then the sum of weights of cones containing $x$ is $\sum_{B\leq I(x)} \mu(\eset,B)$.  
This equals $0$ if $I(x)\neq\eset$, but if $x$ is proper, then the total weight is $\mu(\eset,\eset)=1$.

Thus we must prove that $x$ satisfies the inequalities of the cone $C(\LOOP(\Phi_0/B))$.  
Let $\beta$ be the top-vertex switching function for $B$.  
In a component $A$ of $\Phi|B$ let $v_k$ be any vertex, let $v_j$ be a top vertex, and let $v_i$ be a top vertex of the component of $\Phi|I(x)$ that contains $v_k$ and $v_j$.  
By the definition of a top-vertex switching function, $g_{kj}=\beta_k-\beta_j=\beta_k$.  As $x$ is a 
potential for $I(x) \supseteq B$, $g_{kj}=x_j-x_k$.  Consequently, $x_j = x_k+\beta_k > \beta_k$ (because $x$ is in the positive orthant).  
By the definition of $\Phi_0/B$, $h_j = \max \beta_k$ over all vertices in $A$.  
It follows that $x_j > h_j$; that is, $x$ is in the cone.
 \end{proof}

There is surely a reciprocity theorem analogous to Stanley's for ordinary graphs \cite{AOG}, based on Ehrhart reciprocity (see \cite[Section 4.6]{EC1}), but merely to state such a result would require the theory of orientation of gain graphs developed in Slilaty's thesis \cite{Slilaty}, which is too large a topic to take up here.  We leave this as a research problem.

\section{Interval graph coloring}\mylabel{graph}
\dateadded{2005 Aug 2}

Our theory applies to a kind of generalized graph coloring.  For each vertex $v_i$ of $\G$ we specify an interval $(h_i,\M] = \{h_i+1, \ldots, \M\} \subseteq \bbZ$, depending on a variable $\M$ but with constant lower end $h_i+1$.  The object is to count the number of proper colorations of $\G$ with $v_i$ colored from the interval $(h_i,\M]$, propriety meaning that no two adjacent vertices have the same color.  Let $\chi_{\G,h}(\M)$ denote this number.

We define $\Pi(\G)$ to be the set of partitions $\pi$ of $V$ such that each block of $\pi$ induces a connected subgraph of $\G$.  The total partition of $V$, whose blocks are singletons, is denoted by $\hat0$.

\begin{cor}\mylabel{C:graph}
Let $\G$ be a graph and let $h: V \to \bbZ$.  For $\M \in \bbZ$,
$$
\chi_{\G,h}(m) = \sum_{\pi\in\Pi(\G)} \mu(\hat0,\pi) \prod_{W\in\pi} [\M - h(W) ]^+ ,
$$
where $\mu$ is the M\"obius function of $\Pi(\G)$ and $h(W) := \max \{h_i: v_i \in W\}$.
\end{cor}

\begin{proof}
This is the special case of the general theory in which all the edge gains equal 0.  We write $0\G$ for $\G$ with all $0$ gains and we define $(0\G,h)$ to be the rooted integral gain graph $\Psi$ such that $\Psi \setminus v_0 = 0\G$ and $v_0 v_i$ has gain set $(-\infty,h_i]$.  Then Theorem \myref{L:rootedintegralexpansion} applies.  A closed, balanced edge set $B$ is uniquely determined by its partition $\pi(B)$, and the partitions $\pi(B)$ are precisely the partitions in $\Pi(\G)$; thus $\mu(\eset,B) = \mu(\hat0,\pi)$.  There is no switching, so $h_B(W) = \max \{ h_i: v_i \in W \}$.
\end{proof}

One can think of $(0\G,h)$ as an integrally weighted ordinary graph.  We note the similarity of our contraction formula for $h$ to that of Noble and Welsh \cite{NW}; informally, their theory differs in having $h(W) = \sum \{ h_i: v_i \in W \}$.  We explore the analogy in \cite{PTRI}.

Corollary \myref{C:graph} generalizes to any rooted integral gain graph in which $\Psi \setminus v_0$ is balanced.  It is only necessary to execute top-vertex switching so that $\Psi \setminus v_0$ has all zero gains.  Assuming for simplicity that $\Psi \setminus v_0$ is connected, switching changes $h_j$ to $h_j^\eta = h_j + \eta_j$ where $\eta_j$ is the gain of a path in $\Psi \setminus v_0$ from $v_j$ to the top vertex.  (All $\eta_j$ are nonnegative.)

\section{Modular coloring}\mylabel{modular}

\emph{Modular coloring} means that we interpret the gains modulo $\M$ and take colors in the group $\bbZ_\M$.  Let us write $\chimod_\Phi(\M)$ for the number of ways to do this, the \emph{modular chromatic function}. 
As we saw, the characteristic polynomial of an affinographic arrangement $\cA$ equals the balanced chromatic polynomial $\chib_{\Phi(\cA)}(\l)$.  This polynomial is unchanged if we take the gains in $\bbZ_\M$ for sufficiently large $\M$; indeed, it suffices that $\M>\max \phi(C)$, the largest gain of a circle, because then the list of balanced circles is certain to remain the same.  Thus we can compute $\chib_{\Phi(\cA)}(\M)$, hence $p_\cA(\l)$, by counting colorations of $\Phi$ with color set $\bbZ_\M$ for $\M>\max \phi(C)$.  This is the approach of Athanasiadis, but explained in the language of gain graphs, in which it is not necessary to restrict $\M$ to be a prime power as Athanasiadis did.

It follows that the modular chromatic function is a polynomial for $\M>\max \phi(C)$.  Another clear picture of why that is so is given by a simplification of the method we applied to integral coloring, simply omitting the root vertex.  Let us write
 \begin{equation*}\mylabel{E:modcoloringsum}
 \bL(\Phi) := \sum_{B \in \Latb(\Phi)} \mu(\eset,B) \LOOP(\Phi/B) ,
 \end{equation*}
 where $\mu$ is the M\"obius function of $\Latb(\Phi)$.  (It is not necessary to switch by the top-vertex rule for contraction; any switching function will yield the same result.)  Each of the graphs $\L=\LOOP(\Phi/B)$, unless it has no edges, has loops with integral gains.  When we take gains modulo $\M$, some loops may find themselves with gain $0$; if this happens, then the contribution of $\L$ to $\chimod_\Phi(\M)$ is zero; but if not, then its contribution is $\M^{|V(\L)|}$ since each vertex can have any color in $\bbZ_\M$.  Thus, for instance, $\chimod_\Phi(\M) = \chib_\Phi(\M)$ for all $\M$ greater than the maximum gain of any loop in any $\L$, which is the same as the largest gain of a circle in $\Phi$, and also for any smaller value of $\M$ that does not divide the gain of any loop; and for no other value of $\M>0$.  Summarizing the essential point:

\begin{thm}\mylabel{T:modcoloring}
 Let $\Phi$ be an integral gain graph and $\M>0$.  The number of ways to properly color $\Phi$ with colors in $\bbZ_\M$ is the evaluation of $\bL(\Phi)$ obtained by substituting for each graph $\LOOP(\Phi/B)$ the value $0$ if it has a loop whose gain is a multiple of $\M$ and $\M^{|V(\L)|}$ otherwise.
 \end{thm}

We may regard each graph in $\bL(\Phi)$ as a product of single-vertex graphs (that is, multiplication is disjoint union) and define the evaluation as a ring homomorphism whose value on a single-vertex graph is $0$ if the vertex supports a loop with gain divisible by $\M$, and $\M$ otherwise.

Let us denote by $\Phi_\M$ the gain graph $\Phi$ with gains interpreted modulo $\M$.  Since $\chimod_\Phi(\M) = \chib_{\Phi_\M}(\M)$ and the polynomial $\chib_{\Phi_\M}(\l)$ satisfies the deletion-contraction identity with respect to any link \cite[Corollary III.3.3]{BG}, it follows that 
 $$
 \chimod_\Phi(\M) = \chimod_{\Phi\setminus e}(\M) - \chimod_{\Phi/e}(\M) \quad\text{for } \M=1,2,\ldots .
 $$

\section{Examples}\mylabel{examples}

For $a$ and $b$ integers with $a \leq b$, define $[a,b]{\vec K_n}$ to be the complete graph with, on each edge $e_{ij}$ for $i<j$, all the gains in the interval $[a,b]:=\{a,a+1,\ldots,b\} \subset \bbZ$.  
It is hard to solve these examples with Theorem \myref{T:integralexpansion} because the M\"obius function is difficult; so we employ coloring.  The methods we use are adapted from the ideas of Athanasiadis \cite{Athan}, who colored cyclically.

\subsection{The Shi and extended Shi arrangements}\mylabel{shi}

The Shi arrangement $\cS_n$ corresponds to the gain graph $[0,1]{\vec K_n}$. We count the integral $\M$-colorations for some positive integer $\M$.  
We consider a coloration as a placement of $n$ distinguishable vertices into $\M$ possible positions.  
The $0$ gains on the edges prevent two vertices from having the same position.  Given two vertices $v_i$ and $v_j$ with $i<j$, the edges with gain $1$ prevent $v_i$ from being immediately before $v_j$.  
More generally, if $l$ vertices $v_{i_1},\ldots, v_{i_l}$ appear in consecutive positions, the 
gain $1$ edges imply that they are disposed in reverse order.  This gives a bijection between proper colorations and distributions of the $n$ vertices into the $\M-n+1$ spaces between and around the $\M-n$ free positions, since the set of vertices in any one space must be in descending order.  
It follows that the number of lattice points in $[\M]^n$ but not in any hyperplane of the Shi arrangement is equal to 
 \begin{equation}\mylabel{E:Shi}
 \pZ_{\cS_n}(\M) = \chiZ_{[0,1]{\vec K_n}}(\M) = \begin{cases} 
 (\M-n+1)^n    &\text{if } \M\ge n ,\\ 
 0             & \text{if } \M < n .
 \end{cases}
 \end{equation}

The extended Shi arrangement $\cS_n(s)$ corresponds to the gain graph $[-s+1,s]{\vec K_n}$, whose integral chromatic function we abbreviate as $\chiZ_s(\M)$.  In order to evaluate this function, we prove a general reduction formula:
 \begin{equation}\mylabel{E:shrink}
 \chiZ_{[-a,b]{\vec K_n}}(\M)=\chiZ_{[0,b-a]{\vec K_n}}(\M-[n-1]a)
 \end{equation}
 when $0 \leq a \leq b$.  We again consider a coloration as a placement of $n$ distinguishable vertices into $\M$ possible positions.  The $[-a,a]$ gains on the edges prevent two vertices from having positions less than $a+1$ apart.  This implies that between two vertices there must be at least $a$ free positions.  If we erase $a$ of these open positions, we have $\M-(n-1)a$ colors available to color the vertices according to the rules of $[0,b-a]{\vec K_n}$.  The conversion is reversible; this proves the formula.

In particular, $\chiZ_s(\M) = \chiZ_1(\M-(n-1)(s-1))$.  We conclude that 
the number of lattice points in $[\M]^n$ but not in any 
hyperplane is
 $$
 \pZ_{\cS_n(s)}(\M) = \chiZ_s(\M) = \begin{cases} 
 [\M-s(n-1)]^n   &\text{if } \M\ge n+(s-1)(n-1) ,\\ 
 0               &\text{if } \M < n+(s-1)(n-1) .
 \end{cases}
 $$
 Obviously, this is a piecewise polynomial, in a paltry way.

\subsection{The Linial and related arrangements}\mylabel{linial}

The Linial arrangement is the case $a=b=1$.  We solve it by reduction to $[a,b]=[0,2]$.

All the cases $[0,b]{\vec K_n}$ (with $b\geq0$) work as follows.  Position the $n$ vertices at different positions along a line.  This is the same thing as taking a permutation $\tau$ of the vertices.  Once they are placed, to ensure a proper coloration we must add $b$ empty colors between each increasing pair of consecutive vertices, which correspond to an ascent in $\tau$.  So if there are $r$ ascents in $\tau$, this first placement takes up exactly $n+br$ colors.  The other $\M-(n+br)$ free colors must be placed in the $n+1$ spaces ($n-1$ of which are already partly occupied by empty colors) delimited by the $n$ vertices.  So every permutation with $r$ ascents gives exactly 
$$\binom{\M-(n+br)+n}{n}$$ 
proper colorations if $\M-(n+br)+n$ is nonnegative, and none otherwise.  The number of permutations with $r$ ascents is the Eulerian number $A(n,r+1)$, so when $b\geq1$ we have the formula
 \begin{equation}\mylabel{E:nonneg}
 \pZ(\M) = \chiZ_{[0,b]{\vec K_n}}(\M) = \begin{cases} 
 \sum_{r=0}^{\lfloor \M/b \rfloor} A(n,r+1)\binom{\M-br}{n} &\text{if } 
\M\geq0,\\
 0    &\text{if } \M < 0.
 \end{cases}
 \end{equation}
 for the affinographic arrangement $\{ x_j=x_i+g : i<j,\ g=0,1,\ldots,b \}$.  This becomes a polynomial when $\M \geq b(n-1)$.

Combining Equations \eqref{E:shrink} and \eqref{E:nonneg} leads to the evaluation
 $$
 \chiZ_{[-a,b]{\vec K_n}}(\M) = \begin{cases}
 \sum_{r=0}^{\lfloor (\M-[n-1]a)/b \rfloor} 
A(n,r+1)\binom{\M-[n-1]a-br}{n} &\text{if } \M\geq(n-1)a ,\\
 0    &\text{if } \M < (n-1)a ,
 \end{cases}
 $$
 when $0 \leq a \leq b$.

There is a transformation which gives the relation
 \begin{equation}\mylabel{E:shift}
 \chiZ_{[0,b]{\vec K_n}}(\M)=\chiZ_{[1,b-1]{\vec K_n}}(\M-n+1) \quad\text{ if } m \geq n-1,
 \end{equation}
assuming $b\geq1$.  (This is the same as \eqref{E:shrink} with $a=-1$ and $b\geq1$.) 
 Due to the zero-gain edges, a coloration of the first graph must use different colors for every vertex, but that need not be so in the second graph.  The other difference in the colorations is that, between two 
vertices with colors $c(v_i)<c(v_j)$ and $i<j$, there must be at least $b$ unused colors for the first graph and only $b-1$ colors for the second.  
We transform a coloration of the first graph to one of the second graph by placing the vertices in the color set $[\M]$ according to the color $c(v_i)$, then taking out a color between two consecutive vertices, except that when the colors are successive integers in the first graph they 
become equal in the second.  That is, the $i$-th vertex in the natural order of the color set $[\M]$ is moved to the left by $i-1$ positions, so $n-1$ colors have been deleted.  (That is possible only if $m\geq n-1$.)  This transformation is a bijection of proper colorations; thus we have Equation \eqref{E:shift}.  
(A nice special case is the Shi-arrangement formula \eqref{E:Shi}, which is an $n$-th power because when $b=1$ the right-side graph of \eqref{E:shift} has no edges.)

From Equations \eqref{E:nonneg} and \eqref{E:shift} we deduce that, for the affinographic arrangement $\{ x_j=x_i+g : i<j,\ g=1,\ldots,b-1 \}$ with $n\geq2$ and $b\geq1$,
 $$
 \pZ(\M) = \chiZ_{[1,b-1]{\vec K_n}}(\M) =
  \sum_{r=0}^{\lfloor (\M+n-1)/b \rfloor} A(n,r+1)\binom{\M+n-1-br}{n} \quad\text{ when } \M\geq0.
 $$
This function is a piecewise polynomial in $\M\geq0$ which becomes 
a polynomial when $\M \geq (b-1)(n-1)$.  \dateadded{17 Feb 2006 corrected}
The Linial arrangement being the case $b=2$, it satisfies the formula 
 $$
 \pZ_n(\M) = \chiZ_{1{\vec K_n}}(\M) =
  \sum_{r=0}^{\lfloor (\M+n-1)/2 \rfloor} A(n,r+1)\binom{\M+n-1-2r}{n} 
 $$
 when $\M\geq0$.  
 
 Inspired by a question from a referee, we noticed a linear factor that appears half the time.  The first few polynomials for large values of $\M$ factor over the reals as
 \begin{align*}
 \pZ_1(\M) &= \M ,  \\
 \pZ_2(\M) &= \M^2-\M+1,  \\
 \pZ_3(\M) &= (\M-1) (\M^2-2\M+4) . 
 \end{align*}
The odd-order polynomials have a factor $\M-\frac{n-1}{2}$.  To see why, let $n=2k+1$.  By the symmetry $A(n,r+1)=A(n,n-r)$, the terms come in pairs $A(n,r+1)\big[ \binom{\M+2k-2r}{2k+1} + \binom{\M-2k+2r)}{2k+1} \big]$ for $r=0, 1, \ldots, k$.  Substituting $\M=k$, the members of each pair have equal magnitude and opposite sign; so they cancel and the sum is zero.  
 \dateadded{factorization 17, 18 Feb 2006}
 
By similar reasoning, the polynomial for each $b\geq2$ and odd $n=2k+1$ has the integral zero $(b-1)k$.  We are inclined to doubt the presence of any other real zeros.
 \dateadded{17, 18 Feb 2006}


\end{document}